\documentclass[11pt,reqno]{amsart}
\usepackage{amsmath,amssymb,bbm,times}

\newcommand{\B}{{\mathbb{B}}}
\newcommand{\C}{{\mathbb{C}}}
\newcommand{\capac}{\operatorname{cap}}
\newcommand{\charfn}{\mathbbm{1}}

\newcommand{\De}{{\Delta}}
\newcommand{\dis}{\displaystyle}
\newcommand{\ep}{{\epsilon}}
\newcommand{\ga}{{\gamma}}
\newcommand{\Half}{{\mathcal H}}

\newcommand{\Imag}{\operatorname{Im}}
\newcommand{\lz}{{\log |z|}}
\newcommand{\N}{\mathbb{N}}

\newcommand{\p}{{\partial}}
\newcommand{\R}{{\mathbb{R}}}
\newcommand{\Real}{\operatorname{Re}}

\newcommand{\supp}{\operatorname{supp}}

\newcommand{\ze}{{\zeta}}

\newtheorem{majortheorem}{Theorem}
\newtheorem{theorem}{Theorem}[section]
\newtheorem{lemma}[theorem]{Lemma}
\newtheorem{proposition}[theorem]{Proposition}
\newtheorem{corollary}[theorem]{Corollary}
\newtheorem{conjecture}{Conjecture}

\newtheorem*{claim}{Claim}
\newtheorem*{example}{Example}
\newtheorem*{formula}{Formula}

\begin{document}

\title[Moment inequalities for equilibrium measures]
{Moment inequalities for equilibrium measures in the plane}

\author{A. Baernstein II, R. S. Laugesen, and I. E. Pritsker}

\address{Department of Mathematics, Washington University, St. Louis, MO 63130, U.S.A.}
\email{Al@math.wustl.edu}

\address{Department of Mathematics, University of Illinois, Urbana, IL 61801, U.S.A.}
\email{Laugesen@uiuc.edu}

\address{Department of Mathematics, Oklahoma State University, Stillwater, OK 74078, U.S.A.}
\email{Igor@math.okstate.edu}

\thanks{Research of Igor Pritsker was partially supported by the
National Security Agency (grant H98230-06-1-0055), and by the
Alexander von Humboldt Foundation. }

\date{\today}

\keywords{capacity, equilibrium potential, Green's function,
extremal problem, polynomials, zeros, critical points}

\subjclass[2000]{Primary 31A15. Secondary 31A05, 30C15, 30E05.}

\begin{abstract}
The equilibrium measure of a compact plane set gives the steady
state distribution of charges on the conductor. We show that certain
moments of this equilibrium measure, when taken about the
electrostatic centroid and depending only on the real coordinate,
are extremal for an interval centered at the origin. This has
consequences for means of zeros of polynomials, and for means of
critical points of Green's functions.

We also study moments depending on the distance from the centroid,
such as the electrostatic moment of inertia.
\end{abstract}

\maketitle

\vspace{-12pt}
\begin{center} \textit{Dedicated to our friend Fred
Gehring, on the occasion of his 80th birthday.}
\end{center}

\section{\bf Introduction}
\label{introduction}

Let $K$ be a compact non-polar subset of the complex plane $\C$, and
$\mu_K$ be its equilibrium measure. For functions $\phi : K \to \R$,
the integral
$$\int_K \phi(z) \, d\mu_K(z)$$
is called the \emph{$\phi$ moment of $K$}, or of the probability measure
$\mu_K$. For example, when $\phi(z) = |z|^2$ the $\phi$ moment is
the moment of inertia about the origin.

In this paper we take up some problems involving maximizing or
minimizing $\phi$ moments when $\phi$ satisfies certain conditions.
In all of our results the competing sets $K$ will have the same
logarithmic capacity, which, as a normalization, we take to be $1$.
That is:
$$
\capac (K) = 1.
$$
And we shall usually take the \emph{conformal
centroid} $\int_K z \, d\mu_K(z)$ of $K$ to lie at the origin:
$$
 \int_K z \, d\mu_K(z) = 0.
$$

In all of our results, the extremal $\phi$ moment will be achieved
by a line segment of length $4$. When sets in the class have
conformal centroid at the origin, the interval $L$ defined by
$$L = [-2,2]$$
will be among the extremals.

In our two main results the function $\phi$ will in fact depend only
on the real part of $z$. Here are those results.
\begin{majortheorem} \label{Kreal}
Suppose $K\subset \R$ is compact with capacity $\capac (K) = 1$, and
that its conformal centroid is at the origin. Then for every convex
function $\phi : \R \to\R$, we have
\[
\int_K \phi(\Real z) \, d\mu_K(z) \geq \int_L \phi(\Real z) \,
d\mu_L(z) .
\]

Moreover, if $K\setminus L$ has positive capacity and the restriction of $\phi$ to $L$ is
not a linear function, then strict inequality holds.
\end{majortheorem}

\begin{majortheorem} \label{Kconn}
Suppose $K \subset \C$ is compact and connected with capacity
$\capac (K) = 1$, and that its conformal centroid is at the origin.
Then for every convex function $\phi : \R \to\R$, we have
\[
\int_K \phi(\Real z) \, d\mu_K(z)  \leq \int_L \phi(\Real z) \,
d\mu_L(z).
\]

Moreover, if $K \neq L$ and the restriction of $\phi$ to $L$ is not
a linear function, then strict inequality holds.
\end{majortheorem}

Theorem~\ref{Kreal} says that among all compact sets on a line with
the same conformal centroid and the same capacity, the least spread
out set, as measured by convex integral means, is a single interval.
Contrary-wise, Theorem~\ref{Kconn} says that among all plane
continua with the same conformal centroid and the same capacity, the
single interval is the most spread out.

For a lower estimate applicable in both Theorems~\ref{Kreal} and
\ref{Kconn}, we observe
\[
\int_K \phi(\Real z) \, d\mu_K(z)  \geq \phi(0)
\]
by Jensen's inequality, whenever $K$ is compact with conformal
centroid on the imaginary axis ($\Real \int_K z \, d\mu_K(z) = 0$)
and $\phi$ is convex. Equality is attained whenever $K$ is contained
in the imaginary axis. Theorems~\ref{Kreal} and \ref{Kconn} also
require only that the conformal centroid be purely imaginary, but,
for brevity, we shall continue to assume the conformal centroid is
at the origin.

The proof of Theorem~\ref{Kconn} is modeled on the proof of a
theorem of Baernstein \cite{Ba1}, p.139, about maximizing
integral means in certain classes of univalent functions in the unit
disk. The novelty in the present Theorem~\ref{Kconn} is that instead
of working with symmetric decreasing rearrangements on circles, as
in \cite{Ba1}, one must devise a ``*-function'' appropriate to
``symmetric increasing rearrangements'' of functions defined and
unbounded in all of $\R$. The proof of Theorem~\ref{Kreal} follows
the same general strategy as that of Theorem~\ref{Kconn}, but is
simpler, in that no functions need to be rearranged.

Theorem~\ref{Kreal} is motivated by considerations in number theory;
it will be proved in \S\ref{Kreal_proof}. Some consequences will be
presented in \S\ref{Kreal_applications}. Theorem~\ref{Kconn}, to be
proved in \S\ref{Kconn_proof}, arose in an attempt to prove a
conjecture stated in \S\ref{momentsmodulus}. The conjecture asserts,
when both $\phi$ and $\phi^\prime$ are convex, that
\begin{equation}
\label{c} \int_K \phi(\log |z|) \, d\mu_K(z) \leq \int_L \phi(\log
|z|) \, d\mu_L(z)
\end{equation}
provided $K$ is compact and connected with $\capac (K)=1$ and with
its conformal centroid at the origin, and with the origin belonging
to $K$. If true, the conjecture would prove a conjecture of
Pommerenke \cite{Po} about integral means of univalent functions in
the class $\Sigma_0$, and would also prove a conjecture of Pritsker
and Ruscheweyh \cite{PR} about lower bounds for factors of
polynomials.

We will present an example showing that \eqref{c} is false within
the class of all convex functions. But if we add to $K$ the
assumption, stronger than having the conformal centroid at the
origin, that $K$ is symmetric with respect to the origin, then
\eqref{c} is true for all convex $\phi$. This is
Corollary~\ref{symmetric}. It follows from Theorem~\ref{logmoment},
which restates a result of Laugesen \cite{L} (and which is itself a
consequence of Baernstein's integral means result, Theorem~1 in
\cite{Ba1}).

\section{\bf Potential theoretic preliminaries}
\label{prelims}

For potential theoretic notions we shall mostly follow the approach
in \cite{Ran}. Let $K$ be a compact subset of $\C$, and $K^c$ be the
complement of $K$ on the Riemann sphere $\overline{\C}$. For a
measure $\mu$ compactly supported in $\C$, the energy $I(\mu)$ is
defined to be $I(\mu) = \int_{K\times K} \log|z-\ze| \, d\mu(z) \,
d\mu(\ze)$. If $I(\mu) = -\infty$ for every $\mu$ supported on $K$
then $K$ is said to be \emph{polar}. If $K$ is non-polar, then there
is a unique probability measure $\mu_K$ on $K$, called the
equilibrium measure of $K$, which maximizes $I(\mu)$ over all
probability measures $\mu$ on $K$. Clearly $I(\mu_K)$ is a finite
real number, because $K$ is bounded. The capacity $\capac (K)$ of
$K$ is defined to be $e^{I(\mu_K)}$. For polar $K$, define $\capac
(K) = 0$. A general set $E$ is said to be polar if $\capac (K) = 0$
for every compact $K \subset E$.

For non-polar $K$, denote by $g$ the equilibrium potential of $K$.
Then
$$g(z) = g_K(z) = \int_K \log |z-\ze| \, d\mu_K(\ze) , \quad z\in \C.$$
Put $g(\infty) = +\infty$. Then $g$ is harmonic in $K^c$ except at
$\infty$, where $g(z) = \log |z| + o(1)$. By Frostman's Theorem
(\cite{Ran}, p.59),
\[
g \geq I(\mu_K) = \log \capac (K)
\]
everywhere in $\C$, with equality on $K\setminus E$ for some polar
set $E$. The potential $g$ is related to the Green's function of
$K^c$ with pole at $\infty$ by
$$g(z) = \log \capac (K) + g(z,\infty, K^c), \quad z \in \C.$$
See \cite{Ran}, pp.107, 132.

Set
$$\B(R) = \{z\in \C: |z|< R\},\quad \overline{\B(R)} = \{z\in \C: |z| \leq R\},$$
and let
$$a_n = a_n(K) = \frac{1}{n} \int_K \ze^n \, d\mu_K(\ze), \quad n\geq
1.$$
In particular, $a_1$ is the conformal centroid of $K$.

Suppose that $K\subset \overline{\B(R)}$. In the definition of $g$,
take $|z| > R$ and expand the $\log$ in powers of $|\zeta|\leq R$.
We obtain
\begin{equation}
\label{d} g(z) = \log |z| - \Real \sum_{n=1}^{\infty} {a_n}z^{-n},
\quad |z| > R.
\end{equation}

Next, suppose that $K_1$ and $K_2$ are two non-polar compact subsets of $\C$
with the same capacity and the same conformal centroid. Defining the
potentials $g_j(z)=\int_{K_j} \log|z-\zeta| \, d\mu_{K_j}(\zeta)$,
from \eqref{d} it follows that
\begin{equation}
\label{e} g_1(z) - g_2(z) = - \Real \sum_{n=2}^{\infty} b_n z^{-n},
\end{equation}
where $b_n = a_n(K_1) - a_n(K_2)$ and the series converges for $|z|
> R$ when $\overline{\B(R)}$ contains both $K_1$ and $K_2$. Thus
$g_1(z) - g_2(z) = O(z^{-2})$ at $\infty$, and also, $g_1$ and $g_2$
are bounded on compact subsets of $\C$. It follows that the function
$$
w(x) = \int_{-\infty}^{\infty} [g_1(x+iy) - g_2(x+iy)]\,dy, \quad
x\in \R,
$$
is well defined as an absolutely convergent integral, and is
continuous on $\R$. (To see the continuity, split the defining
integral into two parts: inside and outside the disk $\B(2R)$.
Outside the disk, $g_1-g_2$ is represented by the absolutely and
uniformly convergent series \eqref{e}, and hence the integral is
continuous in $x$. Inside the disk, one can first write down the
definitions of $g_1$ and $g_2$ as potentials, and then use Fubini's
theorem and integrate the logarithmic kernel with respect to
Lebesgue measure on the vertical segment inside the disk. This
eliminates the singularity, and thus this part of the integral is
continuous in $x$ too.)

Note that for each complete vertical line $\Gamma$ not passing
through $0$ we have $\int_{\Gamma} z^{-n}\,dz = 0, \ n\geq 2$. With
\eqref{e}, this implies
\begin{equation}\label{ff} w(x) = 0, \quad |x| \geq R.
\end{equation}

To prove Theorems~\ref{Kreal} and \ref{Kconn}, we shall make use of
the following formula.

\begin{formula}
Suppose $K_1$ and $K_2$ are compact non-polar subsets of $\C$ having
the same capacity and the same conformal centroid, and contained in
$\overline{\B(R)}$.

Then for each $a\geq R$ and each function $\phi\in C^2(\R)$, we have
$$\int_{K_1} \phi(\Real z) \, d\mu_{K_1}(z) - \int_{K_2} \phi(\Real z) \, d\mu_{K_2}(z)
= \frac{1}{2\pi}\int_{-a}^a w(x)\phi^{\prime \prime}(x)\,dx.$$
\end{formula}

\begin{proof}
For $a>R,\,b>R$ let $Q = [-a,a]\times [-b,b]$. In the sense of
distributions, we have $\De g = 2\pi \, \mu_K$ in $\C$, where $K$
denotes $K_1$ or $K_2$ and $g$ denotes $g_1$ or $g_2$. See
\cite{Ran}, Theorem 3.7.4. For $\psi\in C^2(\C)$, Green's formula
gives
\begin{align*}
2\pi \int_K \psi \, d\mu_K
& = \int_Q \psi \De g\,dx\,dy \\
& = \int_Q g\,\Delta \psi\,dx\,dy + \int_{\p Q} \{ \psi \p_n g - g
\p_n \psi \} |dz|
\end{align*}
where $\p_n$ denotes outer normal derivative. Thus,
\begin{align*}
& 2\pi \{\int_{K_1} \psi \, d\mu_{K_1} - \int_{K_2}\psi \,
d\mu_{K_2} \} \\
& = \int_Q (g_1-g_2) \Delta \psi\,\,dx\,dy + \int_{\p Q} \{\psi \p_n
(g_1-g_2) - (g_1-g_2) \p_n \psi\} |dz| .
\end{align*}

Write $z=x+iy$ and take $\psi(z) = \phi(x)$. Then
\begin{align*}
& 2\pi \{\int_{K_1} \phi(x) \, d\mu_{K_1}(z) - \int_{K_2}\phi(x) \,
d\mu_{K_2}(z)\} \\
& = \int_Q (g_1-g_2)(z)\phi^{\prime \prime}(x)\,dx\,dy + \int_{\p Q}
\{\phi \p_n (g_1-g_2) - (g_1-g_2) \p_n \phi\} |dz| .
\end{align*}

Fix $a>R$ and let $b\to \infty$. Since $g_1(z) - g_2(z) = O(z^{-2})$
and $\nabla [g_1(z) - g_2(z)] = O(z^{-3})$ as $z\to \infty$, the
integral over $Q$ tends to $\int_{-a}^a w(x)\phi^{\prime
\prime}(x)\,dx$ and the boundary integrals over the horizontal sides
tend to $0$.

Write $Q=Q(b)$ to show the dependence on $b$ and denote the right
hand vertical boundary side by $\p Q^+(b)$. Then by \eqref{e},
$$\lim_{b\to \infty} \int_{\p Q^+(b)} (g_1-g_2)\p_n \phi\, |dz| = -
\phi^\prime(a) \Real  \int_{\R} \sum_{n=2}^{\infty}
b_n(a+iy)^{-n}\,dy.$$ The last term equals $\phi'(a)w(a)$ which, by
\eqref{ff}, is $0$. Thus,
$$
\lim_{b\to \infty} \int_{\p Q^+(b)} (g_1-g_2)\p_n \phi\, |dz| = 0.
$$
The three other vertical boundary integrals likewise have limit
zero. The formula is proved when $a>R. $ By continuity, the formula
also holds for $a=R$.
\end{proof}

Our proofs of Theorems~\ref{Kreal} and \ref{Kconn} will make use of
the following lemmas.

\begin{lemma} \label{strip}
With $K_1,\,K_2$ as in the Formula, suppose that $S$ is a vertical
strip $-\infty < \ga_1 < \Real z < \ga_2 < \infty$. If $\mu_{K_1}(S)
= 0$ then $w(x)$ is concave on $(\ga_1,\ga_2)$. If $\mu_{K_2}(S) =
0$ then $w(x)$ is convex on $(\ga_1,\ga_2)$.
\end{lemma}

\begin{proof}
Assume $\mu_{K_1}(S) = 0$. Let $\phi$ be a nonnegative $C^2$
function on $\R$ with compact support in $(\ga_1, \ga_2)$. Take
$a\in \R$  so large that $a\geq R$ and $(\ga_1, \ga_2)\subset
(-a,a)$. Then in the Formula, the integral over $K_1$ is zero. Since
$\phi \geq 0$ and $\mu_{K_2} \geq 0$, the Formula implies that
$$ 0\geq \int_{-a}^a w(x)\phi^{\prime \prime}(x)\,dx = \int_{\ga_1}^{\ga_2} w(x)\phi^{\prime \prime}(x)\,dx.$$

Since $w$ is continuous on $\R$, the $1$-dimensional version of
Weyl's Lemma \cite{Ho} or \cite{Do}, p.127  shows that $w$ is concave on $
(\ga_1,\ga_2)$. The proof is similar when $\mu_{K_2}(S) = 0$.
\end{proof}

\begin{lemma} \label{hopf}
Suppose that a function $u$ is subharmonic in the upper half plane
$\Half$, is continuous on $\Half \cup \R$ and satisfies $u(x_0) >
u(z)$ for all $z\in \Half$ and some $x_0\in \R$. Then
$$\liminf_{y\to 0+} \frac{u(x_0) - u(x_0 + iy)}{y} > 0.$$
\end{lemma}

\begin{proof}
Let $D$ be the open half disk $\Half \cap \B(x_0,\epsilon)$. Then
$u(x) \leq u(x_0)$ for all $x\in [x_0 - \ep, x_0 + \ep]$ and $u(z) <
u(x_0)$ for all $z$ in the circular part of $\p D$. Let $v$ solve
the Dirichlet problem in $D$ with boundary values $u$. Then $v$ is
nonconstant on $\p D$, hence nonconstant in $D$. Also, $\sup_{\p D}
v = v(x_0)$, and so by the strong maximum principle, $v(x_0) > v(z)$
for all $z\in D$. By Hopf's Lemma, as stated at the top and bottom
of \cite{GT}, p.34, the lim inf in Lemma~\ref{hopf} is positive for
$v$. Since $u\leq v$ in $D$, the lim inf is also positive for $u$.
\end{proof}

\section{\bf Proof of Theorem~\ref{Kreal}}
\label{Kreal_proof}

Let $K$ be a compact subset of $\R$. We assume
also that $\capac (K) = 1$ and that the conformal centroid of $K$ is
at the origin, $\int_K z \, d\mu_K(z) = 0$.

Recall that $L = [-2,2]$. Then $\capac (L) = 1$ and the conformal
centroid of $L$ is at the origin. We shall apply the considerations
of \S\ref{prelims} with $K_1 = L$ and $K_2 = K$.

Write $G = g_1$ and $g = g_2$ for the respective Green's functions
of $L^c$ and $K^c$ with poles at $\infty$. The function $w(x)$
introduced in \S\ref{prelims} is defined on $\R$ by
$$w(x) = \int_{\R} [G(x+is)- g(x+is)]\,ds, \quad x\in \R.$$
Hence
$$w(x) = 2\int_0^\infty [G(x+is)- g(x+is)]\,ds, \quad x\in \R,$$
by symmetry of $G$ and $g$ in the real axis, recalling that $K,L
\subset \R$.

As observed in \S\ref{prelims}, $w$ is continuous on $\R$ and
satisfies $w(x) = 0$ for $|x| \geq R$, where $R$ is so large that
$K$ and $L$ are contained in $\overline{\B(R)}$.

Let $\phi:\R \to \R$ be convex. The second distributional derivative
of $\phi$ is a non-negative Borel measure on $\R$; call it $\nu$.
Via approximation, one sees that the Formula in $\S2$ generalizes to
\begin{equation} \label{cx}
\int_L \phi(\Real z) \, d\mu_L(z)  - \int_K \phi(\Real z) \,
d\mu_K(z) = \frac{1}{2\pi}\int_{\R}w \,d\nu.
\end{equation}

Thus, to prove the inequality in Theorem~\ref{Kreal},  it suffices
to prove that
$$w(x)\leq 0 , \quad x \in \R.$$
To accomplish this, we solve the Dirichlet problem in $\Half$ with
boundary values $w$ on $\R$ (and boundary value $0$ at infinity),
and call the resulting function $w(z) = w(x+iy)$. Then $w$ is
continuous on $\Half \cup \R$, equals $0$ on the real axis near
infinity, and tends to $0$ as $z\to \infty$ in $\Half$. Moreover, we
will show $w$ has the representation
\begin{equation}
\label{aa} w(z) = 2 \int_y^{\infty} [G(x+is) - g(x+is)]\,ds, \quad z
= x+iy\in \Half  \cup \R.
\end{equation}

To see that this representation is valid, call the function on the
right $\widetilde{w}$. Then, see \eqref{e}, $\widetilde{w}$ is
bounded and continuous on $\Half \cup \R$, and equals $w$ on $\R$
and tends to $0$ as $z\to \infty$ in $\Half$. Further, $g$ and $G$
are harmonic in $\Half$ (since $K,L \subset \R$), from which it
follows that $\widetilde{w}_{xx} = 2(G_y - g_y) $. Also,
differentiation of $\widetilde{w}$ twice with respect to $y$ gives
$\widetilde{w}_{yy} = 2(g_y - G_y)$. Thus $\widetilde{w}$ is
harmonic in $\Half$. By uniqueness of solutions to the Dirichlet
problem in the halfplane, we have $\widetilde{w} = w$.

As just noted, \eqref{aa} gives the identity
$$
w_{xx}(z) = 2[G_y(z) - g_y(z)], \quad z\in \Half ,
$$
and also gives
\begin{equation}
\label{cc} w_y(z) = 2[g(z)-G(z)], \quad z\in \Half .
\end{equation}

Set
$$M = \sup_{\Half  \cup \R} w.$$
Then $M \geq 0$, since $w(x) = 0$ for $|x| \geq R$. Suppose that $M
> 0$. Then by continuity of $w$ and the strong maximum principle
there exists $x_0\in \R$ such that $w(x_0) = M$ and $w(z) < w(x_0)$
for each $z\in \Half$. There are two possible locations for
$x_0$. \vspace{.1in}

\noindent Case 1. $x_0 \in (-\infty, -2] \cup [2,\infty)$. Since
$\mu_L$ is supported on $[-2,2]$, Lemma~\ref{strip} in
\S\ref{prelims} implies that $w$ is concave on each open bounded
subinterval of $(-\infty, -2]$, hence is concave on $(-\infty, -2]$.
Since $w(x) = 0$ for all $x\leq -R$, it follows that $w\leq 0$ on
$(-\infty, -2]$. So if $x_0\in (-\infty, -2]$, then $M=0$.
Similarly, if $2\leq x_0 < \infty$ then $M= 0$. This contradicts our
assumption that $M > 0$, and so Case 1 cannot occur.

\vspace{.1in} \noindent Case 2. $x_0\in (-2,2)$. Since $w$ is
harmonic in $\Half$, Lemma~\ref{hopf} implies that
\begin{equation}
\label{dd} \liminf_{y\to 0+} \frac{w(x_0) - w(x_0+iy)}{y} > 0.
\end{equation}

On the other hand, $w(x_0 + iy)$ is a continuous function of $y$ on
$[0,\infty)$ and is differentiable on $(0,\infty)$. By the mean
value theorem, for each $y>0$ there exists $y^*\in (0,y)$ such that
\begin{align*}
\frac{w(x_0)-w(x_0 + iy)}{y}
& = -w_y(x_0 + iy^*) \\
& = 2[G(x_0+iy^*)-g(x_0 + iy^*)] && \text{by \eqref{cc}} \\
& \leq 2G(x_0+iy^*) && \text{since $g \geq 0$} \\
& \to 0
\end{align*}
as $y\to 0+$, because $x_0 \in L$ and $G=0$ on $L$. This contradicts
$\eqref{dd}$, and so Case 2 cannot occur.  The inequality in
Theorem~\ref{Kreal} is proved.

To prove the strict inequality statement, assume $K\setminus L$ has
positive capacity. Since $g$ is harmonic in $\C\setminus \supp
\mu_K$ and is nonnegative and nonconstant there, we have by the
strong minimum principle that $g>0$ on $\C \setminus \supp \mu_K$.
Recalling that $g=0$ on $K \setminus E$ for some polar set $E$, we
deduce $K\setminus E\subset \supp \mu_K$. Hence, if $\supp \mu_K
\subset [-2,2]$ then $K \setminus E \subset [-2,2] = L$, which
implies $K \setminus L \subset E$ is polar, meaning $K \setminus L$
has capacity zero in contradiction to our assumption. Therefore
$\supp \mu_K \not \subset [-2,2]$, so that some $t>2$ exists with
either $\mu_K\left([t, \infty) \right) > 0$ or $\mu_K\left((-\infty,
-t]\right) > 0$. Say the former holds. Take $\phi(x) = (x-t)^{+}$ in
(\ref{cx}). Then the distributional second derivative of $\phi$ is
the unit point mass at $t$, so that
$$w(t) = -2\pi \int_K (x-t)^{+}\,d\mu_K(x) < 0.$$
Thus, $w$ is nonconstant in the closed upper half plane, and also
$w\leq 0$ as we saw above.

If $w(x) = 0$ for some $x\in L$, then we can rerun the Hopf's lemma
argument in the proof of Case 2 to get a contradiction. So, $w(x) <
0$ at every $x\in L$. If $\phi$ is convex on $\R$ and not linear on
$L$, then the corresponding measure $\nu$ satisfies $\nu\left((-2,2)
\right) > 0$. Formula (\ref{cx}) implies $\int_{L} \phi(\Real
z)\,d\mu_L(z) < \int_K \phi(\Real z)\,d\mu_K(z)$. $\quad \square$

\section{\bf Applications of Theorem~\ref{Kreal}}
\label{Kreal_applications}

This section contains three direct applications of
Theorem~\ref{Kreal}. They are related to the properties of Green's
function and its derivatives, as well as to the asymptotic zero
distribution of polynomials.

\subsection{Pointwise bounds for Green's function and its
derivatives}

Suppose as before that $K\subset\R$ is a compact set, $\capac (K)=1$
and $\int_K x \, d\mu_K(x)=0$, where $\mu_K$ is the equilibrium
measure of $K$. Recall that $g$ denotes Green's function of
$\overline{\C}\setminus K$ and $G$ denotes Green's function of
$\overline{\C}\setminus L$, with poles at $\infty$, where
$L=[-2,2]$. Then the equilibrium measure of $L$ is given by $d\mu_L
= dx/(\pi\sqrt{4-x^2})$, and $G(z)=\log|z+\sqrt{z^2-4}|-\log 2$.

\begin{corollary} \label{pointbound}
Let $x_0\in\R,\ x_0>2$, be fixed. For any set $K$ as above, with
$\max K < x_0$, we have
\begin{align}
\frac{\partial^m g}{\partial x^m}(x_0) & \leq \frac{\partial^m
G}{\partial x^m}(x_0) \qquad \text{when $m \geq 0$ is even,}
\label{2.1} \\
\frac{\partial^m g}{\partial x^m}(x_0) & \geq \frac{\partial^m
G}{\partial x^m}(x_0) \qquad \text{when $m \geq 1$ is odd.}
\label{2.2}
\end{align}
Furthermore, if $z_0 =x_0+iy_0$ and $\max K < x_0 - |y_0|$, then
\begin{equation}
g(z_0) \le G(z_0). \label{2.3}
\end{equation}
Equality holds in \eqref{2.1}--\eqref{2.3} if and only if
$K\setminus L$ has zero capacity.
\end{corollary}

In words, inequality \eqref{2.1} with $m=0$ says that the Green's
function of $K^c$ is smaller at $x_0$ than the Green's function of
$L^c$, which is reasonable since $K$ is more spread out than $L$ and
thus contains points closer to $x_0$.

Clearly, one can consider $x_0<-2$ by symmetry, and make
corresponding adjustments in the above corollary.

\begin{proof}
Recall that $g(z)=\int_K \log|z-s| \, d\mu_K(s)$. Since $\max K <
x_0$, we have that $g\in C^{\infty}$ around $z=x_0$, and
\[
\frac{\partial^m g}{\partial x^m}(x_0) = \int_K (-1)^{m+1}(m-1)!
(x_0-s)^{-m} \, d\mu_K(s) , \quad m \in \N .
\]
Note that the integrand can be extended to a strictly convex
function of $s \in \R$, for odd $m\in\N$. Hence Theorem~\ref{Kreal}
gives \eqref{2.2}. Similarly, the integrand is strictly concave for
even $m\ge 0$, so that we obtain the reversed inequality
\eqref{2.1}.

For $z_0=x_0+iy_0$, we have
\[
g(z_0)=\frac{1}{2} \int_K \log((x_0-s)^2+y_0^2) \, d\mu_K(s).
\]
Thus the integrand is a strictly concave function of $s$ for $s <
x_0 - |y_0|$, and \eqref{2.3} is again a direct consequence of
Theorem~\ref{Kreal}.

For the case of equality, suppose $\capac(K \setminus L)=0$. Then
$\mu_K(K \setminus L)=0$ (since otherwise the restriction of $\mu_K$
to $K \setminus L$ would give a finite energy), and so $\mu_K$ is
supported in $L$. Hence $g$ is harmonic in $\C \setminus L$, so that
$g-G$ is harmonic in $\overline{\C} \setminus L$. Because $g-G$ is
nonnegative on $L$ and equals zero at infinity, the strong maximum
principle implies $g-G \equiv 0$, so that equality holds in
\eqref{2.1}--\eqref{2.3}. On the other hand, if $\capac(K \setminus
L)>0$ then strict inequalities hold in \eqref{2.1}--\eqref{2.3} by
Theorem~\ref{Kreal}.
\end{proof}

\subsection{Means of zeros of polynomials}

This part is inspired by the problem on the smallest limit point for
the arithmetic means of zeros for polynomials with integer
coefficients and positive zeros, considered by Schur \cite{Sch} and
Siegel \cite{Sie}. They gave lower bounds for the arithmetic means
of zeros, which improved the standard arithmetic-geometric means
inequality.

We consider certain extremal polynomials on the real line here. The
number theoretic aspects of the problem for integer polynomials will
be treated in a separate paper.

Let $K\subset\C$ be an arbitrary compact set. It is well known that
for any monic polynomial $P_n$ of degree $n$, we have $\|P_n\|_K \ge
(\capac (K))^n$, where the norm on $K$ is the supremum norm (cf.
\cite{AB}). Thus a sequence of monic polynomials $P_n,\ n\in\N$, is
called \emph{asymptotically extremal} for the set $K$ if
$$
\lim_{n\to\infty} \|P_n\|_K^{1/n} = \capac (K).
$$
This class includes many polynomials orthogonal with respect to
various weights on $K$, and polynomials minimizing various $L^p$
norms; see \cite{AB} and \cite{StT} for numerous examples. Among the
classical families on the real line, we mention Legendre, Chebyshev
and Jacobi polynomials (normalized to be monic). Asymptotically
extremal polynomials have interesting asymptotic zero distributions.
Let $\{\alpha_{k,n}\}_{k=1}^n$ be the zeros of $P_n$. Define the
counting measure for the set $\{\alpha_{k,n}\}_{k=1}^n$ by
$$
\tau_n = \frac{1}{n} \sum_{k=1}^n \delta_{\alpha_{k,n}},
$$
where
$\delta_{\alpha_{k,n}}$ is a unit point mass at $\alpha_{k,n}$. If
$K\subset\R, \ \capac (K)\neq 0$, and the $P_n$ are asymptotically
extremal for $K$, then the $\tau_n$ form a sequence of positive unit
Borel measures that converge in the weak* topology to the
equilibrium measure of $K$; see Theorem~1.7 of \cite[p. 55]{AB}. The
definition of weak* convergence states that
\[
\lim_{n\to\infty} \int_\C f \,d\tau_n = \int_\R f \, d\mu_K
\]
for any continuous function $f$ on $\C$. This enables us to obtain
information on the limiting behavior of means of zeros of $P_n$. In
particular, we have the following result stated for $K$ normalized
by $\capac (K)=1$. (The case of arbitrary capacity can be reduced to
this by a linear change of variable.)

\begin{corollary} \label{arithmetic}
Suppose that $\phi:\C\to\R$ is continuous, and $\phi$ is convex on
$\R$. Assume that $K\subset\R$ is compact, $\capac (K)=1$ and
$\int_K x \, d\mu_K(x)=0$. If $P_n,\ n\in\N$, is a sequence of
asymptotically extremal polynomials for $K$, then we have for the
$\phi$-arithmetic means of their zeros that
\begin{align}
\lim_{n\to\infty} \frac{1}{n} \sum_{k=1}^n \phi(\alpha_{k,n})
& = \int_K \phi(x) \, d\mu_K(x) \label{1.9} \\
& \geq \int_L \phi(x) \, d\mu_L(x) = \int_{-2}^2
\frac{\phi(x)\,dx}{\pi\sqrt{4-x^2}} =: \ell(\phi). \notag
\end{align}
If $K\setminus L$ has positive capacity and the restriction of
$\phi$ to $L$ is not a linear function, then strict inequality
holds.
\end{corollary}

In particular, if $\phi(x)=|x|^m,\ m\in\N$, then
\[
\ell(|x|^m) = 2^m \frac{\Gamma(m/2+1/2)}{\sqrt{\pi}\,\Gamma(m/2+1)}
,
\]
because the change of variable $x=2t^{1/2}$ reduces the integral for
$\ell(|x|^m)$ to a beta integral. Hence $\ell(|x|)=4/\pi$ and
$\ell(x^2)=2$.

\begin{proof}
Since $\phi$ is continuous on $\R$, the first equality in
\eqref{1.9} follows from the weak* convergence of $\tau_n$ to
$\mu_K$. The inequality (and when it becomes equality) is immediate
from Theorem~\ref{Kreal}.
\end{proof}

We also state a version of this result for polynomials with positive
zeros.

\begin{corollary} \label{positivezeros}
Assume $\phi:[0,\infty)\to\R$ and that $\phi(x^2)$ is convex on
$\R$. Suppose $K\subset[0,\infty)$ is compact and $\capac (K)=1$. If
$P_n,\ n\in\N$, is a sequence of asymptotically extremal polynomials
for $K$, and if each $P_n$ has all its zeros positive, then
$$
\lim_{n\to\infty} \frac{1}{n} \sum_{k=1}^n \phi(\alpha_{k,n}) =
\int_K \phi(x) \, d\mu_K(x) \ge \int_0^4
\frac{\phi(x)\,dx}{\pi\sqrt{x(4-x)}} =: \ell_+(\phi).
$$
If $K\setminus [0,4]$ has positive capacity and the restriction of
$\phi(x^2)$ to $[0,4]$ is not a linear function, then strict
inequality holds.
\end{corollary}

In particular, setting $\phi(x)=x^m,\ m\in\N$, gives
$$
\ell_+(x^m) = \int_0^4 \frac{x^m\,dx}{\pi\sqrt{x(4-x)}} = 2^m\,
\frac{1\cdot 3\cdot\ldots\cdot(2m-1)}{m!}.
$$
The first few values of $\ell_+(x^m)$ are $2$ for $m=1$, $6$ for
$m=2$, and $20$ for $m=3$.

\begin{proof}
The proof is essentially the same as for Corollary~\ref{arithmetic}.
For the inequality, one should apply the change of variable $x=t^2$
and define the compact set $\sqrt{K}=\{t\in\R: t^2\in K\}$. Then
$\sqrt{K}$ is symmetric about the origin, so that $\int_{\sqrt{K}} t
\, d\mu_{\sqrt{K}}(t)=0$. Furthermore, $d\mu_{\sqrt{K}}(t) =
d\mu_K(t^2),\ t\in \sqrt{K}$, and $\capac (\sqrt{K})=1$; see
\cite[p. 134]{Ran}. Now apply Theorem~\ref{Kreal} to $\sqrt{K}$.
\end{proof}

A consequence of Corollary~\ref{positivezeros} is that we also have
information on the asymptotic behavior of the coefficients of $P_n$.
For example, if $P_n(x)= x^n + a_{n-1,n} x^{n-1} + \ldots + a_{0,n}
= \prod_{k=1}^n (x-\alpha_{k,n})$ then $a_{n-1,n} = - \sum_{k=1}^n
\alpha_{k,n}$. Hence
\[
\lim_{n\to\infty} \frac{a_{n-1,n}}{n} = - \int_K x \, d\mu_K(x) \le
-2
\]
under the assumptions of Corollary~\ref{positivezeros}, with
equality for $K=[0,4]$.

\subsection{Equilibrium measure and Green's function when
$K$ is the union of several intervals}

Let $K=\bigcup_{l=1}^N\, [a_l,b_l]$, where
$a_1<b_1<a_2<b_2<\ldots<a_N<b_N$ are real numbers. Define the
function $R(z)=\prod_{l=1}^N (z-a_l)(z-b_l)$. Consider the branch of
$\sqrt{R(z)}$, satisfying $\lim_{z\to\infty} \sqrt{R(z)}/z^N = 1$,
which is analytic in $\C\setminus \bigcup_{l=1}^N [a_l,b_l]$. For
future reference, we describe the values of $\sqrt{R(z)}$ on the
real line:
\begin{equation} \label{5.7}
\sqrt{R(x)} = \left\{
\begin{array}{ll}
\sqrt{|R(x)|}, \quad &x\ge b_N, \\
(-1)^{N+l}\, i\, \sqrt{|R(x)|}, \quad &a_l\le x \le b_l,\ l=1,\ldots,N, \\
(-1)^{N+l}\, \sqrt{|R(x)|}, \quad &b_l\le x \le a_{l+1},\ l=1,\ldots,N-1, \\
(-1)^N\, \sqrt{|R(x)|}, \quad &x \le a_1.
\end{array}
\right.
\end{equation}
Here, the values of $\sqrt{R(x)}$ for $x \in \bigcup_{l=1}^N
[a_l,b_l]$ are the limit values of $\sqrt{R(z)}$ when $\Imag z \to
0^+$.

When $K = L = [-2,2]$, then $R(z) = z^2-4$ and for $-2<x<2$ we have
$$d\mu_L(x) = \frac{dx}{\pi\sqrt{4-x^2}} = \frac{dx}{\pi i \sqrt{R(x)}}.$$

We give the following explicit representation for the equilibrium
measure of the set $K$ (see also \cite{StT} and \cite{Tot}).

\begin{proposition} \label{explicitrep}
Let $K=\bigcup_{l=1}^N\, [a_l,b_l] \subset \R$. There exists a
polynomial $T(x)=-x^{N-1}+\ldots\in\R_{N-1}[x]$, such that the
equilibrium measure of $K$ is given by
\begin{equation} \label{5.9}
d\mu_K(x)=\frac{T(x)\,dx}{\pi i \sqrt{R(x)}}, \qquad x \in
\bigcup_{l=1}^N\, [a_l,b_l].
\end{equation}
Furthermore, when $N\geq 2$ we have $T(x)=-\prod_{j=1}^{N-1}
(x-z_j)$ with $z_j\in(b_j,a_{j+1}),\ j=1,\ldots,N-1$, and
\begin{equation} \label{5.10}
\int_K x \, d\mu_K(x) = \sum_{l=1}^N \frac{a_l+b_l}{2} -
\sum_{l=1}^{N-1} z_l.
\end{equation}
\end{proposition}

For the proof of Proposition~\ref{explicitrep}, we need the
following simple lemma.

\begin{lemma} \label{simple}
Let $K=\bigcup_{l=1}^N\, [a_l,b_l]$. For any
$T_{N-1}\in\R_{N-1}[x]$, we have
\begin{equation} \label{5.8}
\frac{1}{\pi i} \int_K \frac{T_{N-1}(t)\,dt}{(t-z) \sqrt{R(t)}} =
\left\{
\begin{array}{ll}
0, \quad &z\in \bigcup_{l=1}^N\, (a_l,b_l), \\
T_{N-1}(z)/\sqrt{R(z)}, \ &z \in \C\setminus K,
\end{array}
\right.
\end{equation}
where the integral is understood in the Cauchy principal value
sense.
\end{lemma}

We remind the reader that when $t\in \R,\;\sqrt{R(t)}$ is defined to
be \newline $\lim_{s\to 0+} \sqrt{R(t+is)}$.

\begin{proof}[Proof of Lemma~\ref{simple}]
For $z\in \C\setminus K$ define $f(z) = T_{N-1}(z)/\sqrt{R(z)}$. It
is easy to see that the limit values of $\sqrt{R(z)}$ as $z$ tends
to a point of $K$ from above and from below are negatives of each
other, so the same is true for $f$. Thus, with obvious notation,
\begin{equation}
\label{4a} f(z+) = f(z) = -f(z-), \quad z\in K.
\end{equation}

Consider a contour $\Gamma$ which consists of $N$ simple closed
curves, one around each of the intervals $[a_l, b_l]$, and located
close to those intervals. Then
$$\frac{1}{2\pi i} \int_{\Gamma} \frac{f(t)}{t-z}\,dt  = f(z)$$
for $z$ in the exterior of $\Gamma$, and for $z\in K$ the integral
equals zero.

Taking $z\in \C\setminus K$, letting $\Gamma$ shrink to $K$, and
using \eqref{4a}, we obtain
$$f(z) = \frac{1}{\pi i}\int_K \frac{f(t)}{t-z}\,dt, \quad z\in K^c,$$
as asserted by the Lemma.

Next, take $z\in \cup_{l=1}^N (a_l,b_l)$. The existence of the
Cauchy principal value at $z$ for the function $f$  follows from the
results in Chapter~2 in [8], which also contains a discussion of
Plemelj's formula. This formula asserts that the Cauchy principal
value satisfies
$$ \frac{1}{\pi i}\int_K \frac{f(t)}{t-z}\,dt = \frac{f(z+) + f(z-)}{2}, \quad z\in \cup_{l=1}^N (a_l,b_l).$$
By \eqref{4a}, the right hand side is zero. This completes the proof
of the lemma.
\end{proof}

\begin{proof}[Proof of Proposition~\ref{explicitrep}]
We shall deduce \eqref{5.9} from Lemma~\ref{simple}. Select
$T(t)=\sum_{j=0}^{N-1} c_j t^j \in \R_{N-1}[t]$ so that it satisfies
the following equations:
\begin{equation} \label{5.11}
\int_{b_l}^{a_{l+1}} \frac{T(t)\,dt}{\sqrt{R(t)}} = \sum_{j=0}^{N-1}
c_j \int_{b_l}^{a_{l+1}} \frac{t^j\,dt}{\sqrt{R(t)}} = 0, \quad
l=1,\ldots,N-1,
\end{equation}
and
\begin{equation} \label{5.12}
\frac{1}{\pi i} \int_K \frac{T(t)\,dt}{\sqrt{R(t)}} =
\sum_{j=0}^{N-1} \frac{c_j}{\pi i} \int_K
\frac{t^j\,dt}{\sqrt{R(t)}} = 1.
\end{equation}
The polynomial $T(t)$ is defined by these equations uniquely,
because the corresponding homogeneous system of linear equations
(with zero on the right of \eqref{5.12}), in the coefficients $c_j$
of $T(t)$, has only the trivial solution. Indeed, let $T_h(t)$ be a
nontrivial solution of this homogeneous system. Since the sign of
$\sqrt{R(t)}$ is constant on each $(b_l,a_{l+1})$ by \eqref{5.7},
$T_h(t)$ must change sign on each $[b_l,a_{l+1}],\ l=1,\ldots,N-1$,
by \eqref{5.11}. Hence $T_h(t)$ has a simple zero in each
$(b_l,a_{l+1}),\ l=1,\ldots,N-1$, and its sign alternates on the
intervals $[a_l,b_l],\ l=1,\ldots,N$. (Note that the same is true
for $T(t)$.) It follows from \eqref{5.7} that $T_h(t)/(\pi i
\sqrt{R(t)})$ doesn't change sign on $K$, contradicting
\[
\frac{1}{\pi i} \int_K \frac{T_h(t)\,dt}{\sqrt{R(t)}} = 0.
\]
Thus $T(t)$ exists and is unique. In addition, the above argument
and \eqref{5.12} show that $T(t)/(\pi i \sqrt{R(t)})$ keeps positive
sign on $K$, that is, \eqref{5.9} actually defines a positive unit
Borel measure on $K$.

As in the proof of Lemma~\ref{simple}, set $f = T/\sqrt{R}$. Let
$$h(z) = \frac{1}{\pi i}\int_K \frac{f(t)}{t-z}\,dt, \quad  z\in \C.$$

Then $h$ is the Cauchy transform of $f\charfn_K$ in $\C$ and is the
Hilbert transform of $f\charfn_K$ on $\R$. It is easy to see that
$f\in L^p(\R)$ for each $1<p<2$. From M. Riesz's conjugate function
theorem (see for example Stein-Weiss \cite{SW}), it follows that
$h\in L^p(\R)$.

From Lemma~\ref{simple}, we see that $h = 0$ on $K$ except at
endpoints, and $h = f$ on $K^c$. Define
$$u(z) = \frac{1}{\pi i}\int_K (\log|z-t|) f(t)\,dt.$$
Then $u$ is continuous on $\C$ and $u_x = -\Real h $ in the open
upper half plane. Since $h\in L^p(\R)$, the function $h(\cdot +iy)$
converges to $h(\cdot)$ in  $L^p(\R)$ when $y\to 0+$, and hence
converge to $h$ in $L^1(a_1, b_N)$.   Thus, for $x\in [a_1,b_N]$,
\begin{align*}
u(x+iy) - u(a_1+iy) & = -\Real  \int_{a_1}^x h(t+iy)\,dt \\
& \to -\Real \int_{a_1}^x h(t)\,dt = -\int_{a_1}^x h(t)\,dt .
\end{align*}

The last equality holds because $h=0$ on $K$ and $h=f$ with $f$ real
in the gaps between the intervals of $K$. Combining this description
of $h$ with \eqref{5.11}, we see that if $x\in K$ then the last
integral is zero. Since $u(x+iy) \to u(x)$ as $y\to 0+$ for all $x$,
we conclude that $u$ is constant on $K$. Then Frostman's theorem and
the uniqueness of the equilibrium measure imply that $f(x)/(\pi
i)\,dx = T(x)\,dx/(\pi i \sqrt{R(x)})$ is the equilibrium measure
for $K$.

We now show that the leading coefficient of $T$ is $-1$. Observe
that \eqref{5.8} gives for $z=0$ and $T_{N-1}(x)=x^{j+1}$ that
\[
\frac{1}{\pi i} \int_K \frac{t^j\,dt}{\sqrt{R(t)}} = 0, \qquad
j=0,\ldots,N-2.
\]
Also, recall that near infinity
\[
\left(1-\frac{a}{z}\right)^{-1/2} = 1 + \frac{1}{2} \frac{a}{z} +
\ldots.
\]
Therefore, we have the following Laurent expansion at infinity
\begin{equation} \label{5.13}
\frac{z^N}{\sqrt{R(z)}} = 1 + \frac{1}{2} \sum_{l=1}^N (a_l+b_l)
\frac{1}{z} + \ldots.
\end{equation}
Applying the same argument as in the proof of Lemma~\ref{simple} and
evaluating the residue at infinity by \eqref{5.13}, we obtain
\[
\frac{1}{\pi i} \int_K \frac{t^{N-1}\,dt}{\sqrt{R(t)}} = -
\frac{1}{2\pi i} \oint_{|z|=r} \frac{z^{N-1}\,dz}{\sqrt{R(z)}} = -1.
\]
Hence \eqref{5.12} gives $c_{N-1} = -1$. Similarly, we have
\begin{align*}
\int_K x \, d\mu_K(x) & = - \frac{1}{\pi i} \int_K
\frac{x^N\,dx}{\sqrt{R(x)}} + c_{N-2} \frac{1}{\pi i} \int_K
\frac{x^{N-1}\,dx}{\sqrt{R(x)}} \\
& = \sum_{l=1}^N \frac{a_l+b_l}{2} - \sum_{l=1}^{N-1} z_l,
\end{align*}
because $c_{N-2}=\sum_{l=1}^{N-1} z_l$.
\end{proof}

 We remark that the zeros of the polynomial $T$ are
exactly the critical points of the Green's function
$g(z,\infty,K^c)$ for the domain $K^c=\overline\C \setminus K$, with
pole at infinity. Indeed, we have for $g(z,\infty,K^c)=\int_K
\log|z-t| \, d\mu_K(t) - \log \capac (K)$ that
\[
g_x(x,\infty,K^c) = \frac{1}{\pi i} \int_K \frac{T(t)\,dt}{(x-t)
\sqrt{R(t)}} = - \frac{T(x)}{\sqrt{R(x)}}, \quad x \in \R\setminus
K,
\]
by \eqref{5.8}. Moreover, $g_y(z,\infty,K^c)$ is zero on $\R
\setminus K$ and is never zero on $\C\setminus \R$.

Thus we can obtain interesting information about location of the
critical points. For the background material on the critical points
of Green's function see Chapter VII of Walsh \cite{Wal}. If
$K=[a_1,b_1]\cup[a_2,b_2]$ and $|b_1-a_1|=|b_2-a_2|$, then it
follows by an elementary symmetry argument that $z_1=(b_1+a_2)/2$.
Also, if $|b_1-a_1|>|b_2-a_2|$ then $z_1>(b_1+a_2)/2$. But the
location of critical points becomes difficult to predict for three
or more intervals.

The following inequality gives information on the average position
of the critical points in terms of the midpoints of the gaps between
the intervals of $K$.

\begin{corollary} \label{average}
Let $K=\bigcup_{l=1}^N\, [a_l,b_l] \subset \R$ satisfy $\capac
(K)=1$. With the above notation, we have
\[
\sum_{l=1}^{N-1} \left(\frac{b_l+a_{l+1}}{2} - z_l\right) \ge 2 -
\frac{b_N-a_1}{2} ,
\]
where the sum is interpreted to be $0$ for $N=1$. Equality holds
above if and only if $K$ is a segment of length 4.
\end{corollary}

\begin{proof}
We consider the integral $\int_K x \, d\mu_K(x)$, and observe that
translating the set $K$ by a constant $c\in\R$ changes the integral
by adding $c$. Hence we may assume that $a_1=0$, and must show that
$\int_K x \, d\mu_K(x) \ge 2$, with equality only for $K=[0,4]$.
Define $\sqrt{K}=\{t\in\R: t^2\in K\}$. Then, as in the proof of
Corollary~\ref{positivezeros}, $\sqrt{K}$ is symmetric about the
origin, $\int_{\sqrt{K}} x \, d\mu_{\sqrt{K}}(x)=0$, and $\capac
(\sqrt{K})=1$. Moreover,
\[
\int_{\sqrt{K}} t^2 \, d\mu_{\sqrt{K}}(t) = \int_K x \, d\mu_K(x) =
\sum_{l=1}^N \frac{a_l+b_l}{2} - \sum_{l=1}^{N-1} z_l
\]
by \eqref{5.10}. Applying Theorem~\ref{Kreal} with $\phi(t)=t^2$, we
obtain that
\[
\int_{\sqrt{K}} t^2 \, d\mu_{\sqrt{K}}(t) \ge \int_{-2}^2
\frac{t^2\,dt}{\pi\sqrt{4-t^2}} = 2,
\]
with equality possible only if $\sqrt{K}=[-2,2]$ and $K=[0,4]$.
\end{proof}

Using higher moments will give more complicated inequalities
involving the endpoints of $K$ and zeros (or coefficients) of $T$.

\section{\bf Proof of Theorem~\ref{Kconn}}
\label{Kconn_proof}

Let $K$ be a compact connected subset of $\C$ with $\capac (K) = 1$.
The connectivity of $K$ implies that each boundary point of the
domain $K^c$ is regular for the Dirichlet problem in $K^c$, which,
in turn, implies that the Green function of $K^c$ is continuous in
$\C$.

We shall assume also that the conformal centroid of $K$ is at the origin.
That is:
$$\int_K z \, d\mu_K(z) =  0.$$
Then by Theorem~1.4 of \cite[p.19]{Pobook},  we have
$$K\subset \overline{\B(2)}.$$

Recall that $L = [-2,2]$. Then $\capac (L) = 1$ and the conformal
centroid of $L$ is at the origin. We shall apply the considerations
of \S\ref{prelims} with $K_1 = K,\;K_2 = L$ and $R = 2$.

Write $g = g_1$ and $G = g_2$ for the respective Green's functions
of $K^c$ and $L^c$ with poles at $\infty$. The function $w(x)$ is
defined on $\R$ by
\begin{equation}
w(x) = \int_{\R} [g(x+is) - G(x+is)]\,ds , \quad x\in \R.
\label{wdef2}
\end{equation}
By the Formula in \S\ref{prelims}, to prove Theorem~\ref{Kconn} it
suffices to prove that
$$w(x)\leq 0 , \quad x \in \R .$$
To accomplish this, we shall extend $w$ to a certain function $w(z)$
which is subharmonic in the upper half plane $\Half$.

For sets $E\subset \R$, let $E^b$ denote the complement of $E$ in
$\R$:
$$E^b = \R\setminus E.$$

Also, let $|E|$ denote the one-dimensional Lebesgue measure of $E$,
and for $y\geq 0$, let
$$I(y) = [-y,y].$$

For bounded $E\subset \R$ with $|E|=2y$ and $x\in \R$, set
\begin{align*}
w(x,E)
& = \dis \int_{\R} [\charfn_{E^b}(s) g(x+is) - \charfn_{I(y)^b}(s)G(x+is)]\,ds \\
& = w(x) + \dis \int_{I(y)} G(x+is)\,ds - \dis \int_E g(x+is)\,ds,
\end{align*}
where $\charfn$ denotes a characteristic or indicator function. The
asymptotic behavior of $g$ and $G$ (discussed in $\S\ref{prelims}$)
ensures that the first integral is absolutely convergent. The second
equality follows from \eqref{wdef2}.

Now take $z= x+iy\in \Half$, and define
\begin{equation} \label{wdef1}
w(z)= \sup_E w(x,E) ,
\end{equation}
where the sup is taken over all bounded measurable $E\subset \R$
with $|E| = 2y$.

For each $x$, we have $g(x+is)\geq 0$ and $\lim_{|s| \to\infty}
g(x+is) = \infty$. The analysis on p.149 of \cite{Ba1} is applicable
to $-g(x+is)$ as a function of $s$, and shows that for each $y \in
[0,\infty)$ there exists a set $E\subset \R$ with $|E| = 2y$ for
which the supremum of  $- \int_E g(x+is)\,ds$ over all $E$ with $|E|
= 2y$ is attained. Note the minus sign in $-g$. Moreover, there
exists a number $t\geq 0$  such that $\{s\in \R: g(x+is) < t\}
\subset E \subset \{g(x+is) \leq t\}$, and $E$ is bounded. We shall
denote such a maximizing set by $E(z)$.  Then

\begin{align}
w(z) & = \dis \int_{\R} [\charfn_{E(z)^b}(s) g(x+is) -
\charfn_{I(y)^b}(s)G(x+is)]\,ds \label{5aa} \\
& = w(x) + \dis \int_{I(y)} G(x+is)\,ds - \dis \int_{E(z)}
g(x+is)\,ds \label{5a}\\
& = \int_{|s|>y} [g(x+is) - G(x+is)]\,dx \label{5aaa} \\
& \qquad \qquad + \int_{I(y)} g(x+is)\,ds - \int_{E(z)} g(x+is) \,ds
. \notag
\end{align}

The following lemma provides information on the maximizing sets
$E(z)$.

\begin{lemma} \label{max}
With the situation as above, there exist positive constants $b$ and
$k$ depending only on $K$ such that whenever $z=x+iy \in \Half$:

\noindent (a) if $y \geq b$  then $E(z)=I(y)+t=[-y+t, y+t]$, for
some $t$ with $|t| < k/y$;

\noindent(b) if $y \leq b$ then $E(z)\subset [-2b, 2b]$.
\end{lemma}

\begin{proof}[Proof of Lemma~\ref{max}] By \eqref{d}, we can write
\begin{equation} \label{gdecomp}
g(z) = \log |z| + h(z)
\end{equation}
where $h$ is harmonic outside $\overline{\B(2)}$. The conformal
centroid of $K$ is at the origin, and so the coefficient $a_1$ in
\eqref{d} equals $0$. Thus, $h$ satisfies $|h(z)| \leq
\frac{k}{8}|z|^{-2}$ and $|\nabla h(z)| \leq \frac{k}{8}|z|^{-3}$
for $|z| \geq 3$, for some positive constant $k$.

From $g_y = y|z|^{-2}+ h_y$, it easily follows that there exists
$b_0 \geq 3$ such that  $g_y(x+iy) > 0$ whenever $y \geq b_0$ and
$g_y(x+iy) < 0$ whenever $y\leq -b_0$.

Now we establish two estimates:
\begin{align}
|g(z) - g(\overline{z})| & \leq \frac{k}{4} |z|^{-2} , \label{est1} \\
g(z+it) - g(z) & \geq \frac{1}{2} ty|z|^{-2} - \frac{k}{4} |z|^{-2}
, \label{est2}
\end{align}
when $z = x+iy\in \Half, |z| \geq 3$ and $t \in (0, |z|]$. The first
estimate is obvious from \eqref{gdecomp}. The second follows
similarly, because $\Real (it/z) \in (0,1]$ and so
\[
\log |(z+it)/z| \geq \log (1+ \Real (it/z)) \geq \frac{1}{2} \Real
(it/z) = \frac{1}{2} ty|z|^{-2} .
\]

Moreover, there exists a number $b> b_0 \geq 3$ such that $g(x+is)
> g(x+is_0)$ whenever $|s|\geq b$ and $|s_0|\leq b_0$, as we now show. For $|x| \leq
3$ one just takes $b$ large enough that $\max_{S_0} g <
\min_{|x|\leq 3} g(x \pm ib)$ where $S_0 = [-3,3] \times [-b_0,
b_0]$, recalling here that $g$ is continuous and finite in the
plane. For $|x| > 3$, one estimates $|g(x+is) - \log|x+is|| \leq
\frac{k}{8} x^{-2}$ and uses concavity of the function $t \mapsto
\log \frac{1+b^2 t}{1+b_0^2 t}, t \in [0,\infty)$, together with
monotonicity properties of $g$; note that for our purposes,
$t=x^{-2} \in (0,1/9)$. Details are left to the interested reader.

Now fix $x\in \R$ and visualize the graph of $p(s)= g(x+is)$. The
function $p$ is strictly increasing on $[b_0, \infty)$, strictly
decreasing on $(-\infty, -b_0]$, and $p(s) > p(s_0)$ for every
$|s_0| \leq b_0$ and $|s| \geq b$. For $\alpha > 0$, write
$E_{\alpha} = \{s: p(s) < \alpha\}$. Then $E_{\alpha}$ is a maximal
set of measure $|E_{\alpha}|$. Set $\alpha_0 = \min \{p(-b), p(b)\}$
and $y_0 = \frac{1}{2} |E_{\alpha_0}|$. Then $E_{\alpha_0}$ is a
single interval which contains $[-b_0, b_0]$, and $y_0 \leq b$.
Given $y\geq b$, there is a unique $\alpha \geq \alpha_0$ such that
$|E_{\alpha}| = 2y$. Then $E(x+iy) = E_{\alpha}$, and this
$E_{\alpha}$ also is a single interval containing $[-b_0, b_0]$.
These facts imply that $E(x+iy)$ has the form $[-y+t, y+t]$, where
$|t| \leq y-b_0$. Further, the maximality of $E(x+iy)$ and
continuity of $g$ imply that $p(-y+t) = p(y+t) = \alpha$.

Take $z = x + iy \in \Half$ with $y \geq b$. Suppose the number $t$
in the previous paragraph is nonnegative; the case $t\leq 0$ is
handled analogously. Let $z_2 = z+it, z_1 = \overline{z}+it$. Then
\begin{align*}
\frac{1}{2} ty|z|^{-2} - \frac{k}{4} |z|^{-2} & \leq g(z_2) - g(z)
&& \text{by \eqref{est2}} \\
& \leq  g(z_2) - g(\overline{z}) + \frac{k}{4} |z|^{-2} && \text{by \eqref{est1}} \\
& < g(z_2) - g(z_1) + \frac{k}{4} |z|^{-2} && \text{since $g(z_1) =
\alpha
< g(\overline{z})$} \\
& = \frac{k}{4}|z|^{-2}
\end{align*}
because $g(z_2) = g(z_1) = \alpha$ as above. Hence $ty < k$, proving
part (a).

To prove (b), take $z = x + iy \in \Half$ with $y \leq b$, and let
$E(z)$ be a maximizing set for $z$, so that $|E(z)|=2y$. Suppose
$E(z)$ intersects the interval $(2b,\infty)$ in a set of measure
$\ep>0$. Then the set $[0,2b] \setminus E(z)$ has measure at least
$\ep$. Since $g(x+is_0) < g(x+is)$ when $0<s_0< 2b < s$, we can
strictly decrease $\int_{E(z)} g(x+is)\,ds$ if we move $E(z)\cap
(2b,\infty)$ into some subset of $[0,2b] \setminus E(z)$. This
violates the definition of maximizing set, and shows that $E(z)$
cannot intersect the interval $(2b,\infty)$ in a set of positive
measure. Similarly it cannot intersect $(-\infty,-2b)$. Thus $E(z)
\subset [-2b,2b]$, after possibly deleting a set of zero measure
from $E$.

\end{proof}

Here now is the main ingredient in the proof of the theorem.

\begin{claim}
$w$ is subharmonic in $\Half$.
\end{claim}

Let us carry out the proof of Theorem~\ref{Kconn} assuming the
claim.

Firstly, the function $w$ is continuous on $\Half \cup \R$. It is
continuous also at infinity, because $w(z) \to 0$ as $z\to \infty$
in $\Half$, as we now show. From (\ref{5aaa}) it suffices to show
that
$$\lim_{z\to \infty} \int_{|s| > y} [g(x+is) - G(x+is)]\,ds = 0$$
and
$$\lim _{z\to \infty} \left( \int_{I(y)} g(x+is)\,ds - \int_{E(z)} g(x+is)\,ds \right) = 0.$$
The first is a simple consequence of \eqref{e}. The second follows
from Lemma~\ref{max}: when $y \geq b$ use part (a) of the lemma, and
then decomposition \eqref{gdecomp}, and when $y \leq b$  with $|x|
\to \infty$, use part (b) of the lemma and then decomposition
\eqref{gdecomp}. A key fact for the latter case is that $|E(y)| =
|I(y)|$. Details are left to the reader.

Continuing now with the proof of Theorem 2, set
$$M = \sup_{\Half \cup \R} w ,$$
where the supremum is finite because $w$ is bounded at infinity.
Note $M \geq 0$, since from \S\ref{prelims} we know $w(x) = 0$ for
$|x| \geq 2$.

If $M>0$, then by continuity of $w$ and the strong maximum
principle, there exists $x_0\in \R$ such that $w(x_0) = M$ and
$w(z)< w(x_0)$ for each $z\in \Half$. Since $K$ is connected, its
orthogonal projection onto the real axis is a single interval
$[c_1,c_2]$, and since $K\subset \overline{\B(2)}$ we have $[c_1 ,
c_2] \subset [-2,2]$. By Lemma~\ref{strip}, $w(x)$ is concave on
every bounded subinterval of $(-\infty, c_1)$, hence is concave on
$(-\infty, c_1)$. Similarly, $w$ is concave on $(c_2, \infty)$.
Since $w(x) = 0$ for $|x|\geq 2$, we must have $w\leq 0$ on
$\R\setminus [c_1,c_2]$. Thus $x_0\in [c_1,c_2]$. Since $x_0$ is a
maximizing point for $w$, Lemma~\ref{hopf} implies we must have
\begin{equation}
\label{h} \liminf_{y\to 0+}\frac{w(x_0) - w(x_0+iy)}{y} > 0.
\end{equation}

On the other hand, from \eqref{5a}, we see that for $y>0$,
\begin{equation} \label{i}
w(x_0) - w(x_0+iy) = \int_{E(x_0 + iy)} g(x_0+is)\,ds - \int_{I(y)}
G(x_0+is)\,ds.
\end{equation}
Since $E(x_0 + iy)$ maximizes integrals of $-g$, for each $y>0$ and
for each bounded $E\subset \R$ with $|E| = 2y$ we have
\begin{equation} \label{ii}
0 \leq \int_{E(x_0+iy)}g(x+is)\,ds \leq \int_E g(x+is)\,ds.
\end{equation}
Further, because $x_0 \in [c_1,c_2]$ there exists $s_0\in \R$ with
$x_0 + is_0 \in K$, so that $g(x_0 + is_0) = 0$. Taking $E = [s_0 -
y, s_0 + y]$ and using continuity of $g$, we see from \eqref{ii}
that
$$\lim_{y\to 0+} \frac{1}{y} \int_{E(x_0+iy)} g(x_0+is)\,ds = 0.$$
Similarly, $G(x_0) = 0$, and hence $\lim_{y\to 0+} \frac{1}{y}
\int_{I(y)} G(x_0+is)\,ds = 0$. Thus, by \eqref{i},
$$\lim_{y\to 0+} \frac{w(x_0) - w(x_0 + iy)}{y} = 0,$$
which contradicts \eqref{h}.

We conclude that $M>0$ is impossible, and so $M=0$, meaning $w \leq
0$ in $\Half \cup \R$. This completes the proof of the inequality in
Theorem~\ref{Kconn}, modulo the Claim.

To prove the strict inequality statement in the theorem, let $K$ be
a compact set satisfying the hypotheses of Theorem~\ref{Kconn} which
does not coincide with $L$. Then $K$, which is contained in
$\overline{\B(2)}$, cannot contain the points $-2$ or $2$, because
if it did then it would equal $[-2,2]=L$ by the equality case of
\cite[Theorem~1.4]{Pobook}. Hence $-2 < c_1 \leq c_2 < 2$. The
argument that gave $w(t) < 0$ for some $t>2$ in the proof of
Theorem~\ref{Kreal} works again here, except with $K$ and $L$
interchanged, producing that $w(t) < 0$ for every $t \in (-2,2)
\setminus [c_1,c_2]$.

Now we show $w(t) < 0$ for every $t \in [c_1,c_2]$. Suppose instead
that $w(x_0)=0$ for some $x_0 \in [c_1,c_2]$. Note $w \leq 0$ is not
identically zero in $\Half$, by the preceding paragraph, and so
$w<0$ in $\Half$ by the strong maximum principle. Now rerun the
argument used above to rule out the case $M>0$, to obtain a
contradiction. Hence $w<0$ on $[c_1,c_2]$.

We have shown $w<0$ on $(-2,2)$, and so formula (\ref{cx}) (with $K$
and $L$ interchanged) implies the strict inequality that
\[
\int_K \phi (\Real z) \,d\mu_K(z) < \int_L \phi(\Real z)\,d\mu_L(z)
,
\]
when $\phi$ is convex on $\R$ and is not a linear function on
$[-2,2]$.

\begin{proof}[Proof of the Claim] Fix $z = x + iy\in
\Half$. Let $E(z)$ be a corresponding maximal set of measure $2y$,
as in \eqref{5a}. For brevity, we'll write
$$E(z) = E$$
and also
$$I(y) = I.$$

Then \eqref{5aa} says
\begin{equation}
\label{j} w(z) = \int_{\R} [\charfn_{E^b}(s) g(x+is) -
\charfn_{I^b}(s) G(x+is)]\,ds.
\end{equation}

Take $\rho \in (0,y)$. To prove subharmonicity of $w$ it suffices to
show that $w(z)$ is less than or equal to the mean value of $w$ over
the circle with center $z$ and radius $\rho$.

The function $g$ is subharmonic in $\C$ and $G$ is harmonic in
$\C\setminus [-2,2]$. Thus
$$g(x + is)\leq \frac{1}{2\pi} \int_0^{\pi} [g(x+is+\rho e^{i\phi})+ g(x+is+\rho e^{-i\phi})]\,d\phi, \quad s\in \R.$$
If $|s| > y$, then equality holds when $g$ is replaced by $G$.
Substitute the inequality and equality into \eqref{j}, and switch
the order of integration on the right. This gives the inequality
\begin{equation}
\label{k} 2\pi w(z)\leq \int_0^{\pi} [J(\phi) + J(-\phi)] d\phi,
\end{equation}
where
$$J(\phi) = \int_{\R} [\charfn_{E^b}(s)g(x+is+\rho e^{i\phi})-\charfn_{I^b}(s)G(x+is+\rho e^{i\phi})]
\,ds.$$

Fix $\phi\in [0,\pi]$ and set $\ep = \rho \sin\, \phi$. In
$J(\phi)$, substitute
$$ x + is + \rho e^{i\phi} = x + \rho \cos\,\phi + i(s + \ep),$$
then make the change of variable $t = s+\ep$, and integrate over
$\R$. We obtain
$$J(\phi) = \int_{\R} [\charfn_{E^b + \ep}(t)g(x+\rho \cos\,\phi + it) - \charfn_{I^b + \ep}(t)G(x+\rho \cos\,\phi + it)]\,dt.$$

The same equation holds when $\phi$ is changed to $-\phi$ and $\ep$
to $-\ep$. It follows that, for $\phi\in [0,\pi]$, \
\begin{align*}
J(\phi) + J(-\phi) & = \int_{\R}\{[\charfn_{E^b+\ep} +
\charfn_{E^b-\ep}](t)g(x+\rho \cos\,\phi + it) \\
& \qquad - [\charfn_{I^b+\ep} + \charfn_{I^b-\ep}](t)G(x+\rho
\cos\,\phi + it)\}\,dt.
\end{align*}

The argument on the top half of p.148 of \cite{Ba1} shows that for
our set $E = E(x+iy)$ and for $0 < \ep < y$ there exist bounded
measurable sets $A$ and $B$ in $\R$ such that $|A| = 2(y+\ep),\;|B|
= 2(y-\ep)$ and
$$\charfn_{E+\ep} + \charfn_{E-\ep} = \charfn_A + \charfn_B.$$

Using $\charfn_A = 1-\charfn_{A^b}$, etc., one sees that this
equation also holds when the four sets are replaced by their
complements in $\R$. Furthermore, $(E\pm\ep)^b = E^b \pm \ep$, and,
recalling that $I = I(y)$, one can check directly that
$\charfn_{I^b+\ep} + \charfn_{I^b-\ep} = \charfn_{I(y+\ep)^b} +
\charfn_{I(y-\ep)^b}$. Thus,
\begin{align*}
& J(\phi)+ J(-\phi) \\
& = \int_{\R}\{[\charfn_{A^b} + \charfn_{B^b}](t)g(x+\rho \cos\,\phi
+ it) \\ &\hskip1in- [\charfn_{I(y+\ep)^b} +
\charfn_{I(y-\ep)^b]}(t)G(x+\rho \cos\,\phi + it)\}\,dt \\ &=
\int_{\R} [\charfn_{A^b}(t)g(x+\rho\cos\,\phi + it) -
\charfn_{I(y+\ep)^b}(t)G(x+\rho\cos\,\phi + it)]\,dt \\ &+ \int_{\R}
[\charfn_{B^b}(t)g(x+\rho\cos\,\phi + it) -
\charfn_{I(y-\ep)^b}(t)G(x+\rho\cos\,\phi + it)]\,dt \\ &\leq
w(x+\rho\cos\,\phi + i(y+\ep)) + w(x+\rho\cos\,\phi + i(y-\ep)) ,
\end{align*}
by the definition of $w$ as a supremum, in \eqref{wdef1}.
Substitution in \eqref{k} gives
\begin{align*}
2\pi w(z) & \leq \int_0^{\pi} [w(x+\rho\cos\,\phi + i(y+\ep)) +
w(x+\rho\cos\,\phi + i(y-\ep))]\,d\phi \\
& = \int_0^{\pi} [w(z+\rho e^{i\phi}) + w(z + \rho
e^{-i\phi})]\,d\phi,
\end{align*}
recalling $\ep = \rho \sin\, \phi$. Thus, $w$ satisfies the sub-mean
value property at $z$, and the Claim is proved.
\end{proof}

\section{\bf Moments involving $|z|$}
\label{momentsmodulus}

In Theorem~\ref{Kconn} we obtained sharp upper bounds for moments of
the form $\int_K \phi(\Real  z) \, d\mu_K(z)$, where $K$ is a
continuum satisfying certain hypotheses. In this section, we again
take $K$ to be a continuum, and seek sharp upper bounds for moments
of the form $\int_K \phi(|z|) \, d\mu_K(z)$. It turns out to be
convenient to state the results in terms of $\phi(\lz)$ instead of
$\phi(|z|)$.

Let $K$ be a compact, connected subset of $\C$ that contains the
origin, and satisfies $\capac (K) = 1$. As before, let $K^c =
\overline{\C} \setminus K$. Also, set $\Omega = \{z\in \C: 1/z \in
K^c\}$. Then the plane domain $\Omega$ is the image of the unit disk
$\B(1)$ under a function $f$ belonging to the class $S$ of univalent
analytic functions (conformal mappings) in the disk with $f(0) = 0$
and $f'(0) = 1$. See, for example, \cite{Ran}. The Koebe one-quarter
theorem asserts that $\Omega = f(\B(1))$ contains the disk
$\B(1/4)$, so that
$$K \subset \overline{\B(4)}.$$
The interval
$$\widetilde{L} = [0,4]$$
satisfies all our assumptions on $K$, and shows that the ``4'' on
the previous line is the smallest possible constant.

As in \S\ref{prelims}, let $g(z) = \int_{K} \log|z-\zeta| \,
d\mu_K(\zeta)$ be the Green's function of $K^c$ with pole at
$\infty$. From \eqref{d}, we have
$$g(z) = \lz - \Real \sum_{n=1}^{\infty} a_nz^{-n},\quad |z| > 4,$$
where $a_n = n^{-1}\int_K {\zeta}^n \, d\mu_K$. From the expansion,
it follows that
\begin{align}
\frac{1}{2\pi} \int_{-\pi}^{\pi} g(re^{i\theta})\,d\theta & =
\log r, && r\geq 4, \label{6a} \\
\frac{1}{2\pi} \int_{-\pi}^{\pi} g_r(re^{i\theta})\,d\theta & =
r^{-1}, && r\geq 4 . \label{6b}
\end{align}

Next, we give a representation of logarithmic moments of $\mu_K$ in
terms of integrals involving $g$.

\begin{proposition} \label{representation}
Assume that $K\subset \C$ is compact and connected, with $\capac (K)
= 1$. Suppose $\phi \in C^2(\R)$ is constant near $-\infty$. Then
for each $R\geq 4$ we have
\begin{align*}
\int_K \phi(\lz) \, d\mu_K(z)
& = \frac{1}{2\pi} \int_{\B(R)}g(z) \phi^{\prime \prime}(\lz)|z|^{-2}\,dx\,dy \\
& \qquad \qquad + \phi(\log R) - \phi^\prime (\log R)\log R.
\end{align*}
\end{proposition}

\begin{proof}
As in the proof of the Formula in \S\ref{prelims}, we start with
$2\pi \mu_K = \Delta g$, then apply Green's formula to
the integral on the left, this time in the disk $\B(R)$. Set
$\psi(z) = \phi(\lz)$. Then $\psi$ is constant on circles, and from
\eqref{6a} and \eqref{6b} the boundary terms have the form stated.
Also
$$\Delta \psi(z) = \phi^{\prime \prime}(\lz)|z|^{-2},$$
so the integral over $\B(R)$ has the form stated.
\end{proof}

Our Theorem~\ref{Kconn} takes $K$ to be conformally centered. Our
next theorem drops that assumption, assuming instead that $K$
contains the origin and proving that the logarithmic moments are
maximal when $K$ equals the segment $\widetilde{L}=[0,4]$ with one
endpoint at the origin (rather than $L=[-2,2]$, which is centered at
the origin).

Let $\widetilde{G}$ denote the Green's function of
${\widetilde{L}}^c$ with pole at $\infty$.

\begin{theorem} \label{logmoment}
Suppose $K\subset \C$ is compact, connected, contains the origin,
and has $\capac (K) = 1$. Then for every convex function $\phi:\R
\to \R$, we have
$$\int_K \phi(\lz) \, d\mu_K(z) \leq \int_{\widetilde{L}}\phi(\lz) \, d\mu_{\widetilde{L}}(z).$$
\end{theorem}
This result is due to Laugesen \cite[Corollary~6]{L}. We give below
a brief version of that proof, relying on Baernstein's result on
integral means.
\begin{proof}
When $\phi$ is linear, the theorem holds with equality because
\[
\int_K \lz \, d\mu_K(z) = g(0) = 0
\]
(and similarly for $L$), using that $0 \in K$ by hypothesis and that
every point of $K$ is regular for the Dirichlet problem in $K^c$.
For general convex $\phi$, we can reduce by approximation to the
case where $\phi$ is linear near $-\infty$, and hence to the case
where $\phi \equiv 0$ near $-\infty$. Then by mollification we may
further assume $\phi$ is smooth. Then by
Proposition~\ref{representation}, to prove Theorem~\ref{logmoment}
it suffices to show that for every $r\in (0,\infty)$,
$$
\int_{-\pi}^{\pi} g(re^{i\theta})\,d\theta \leq \int_{-\pi}^{\pi}
\widetilde{G}(re^{i\theta})\,d\theta.
$$

As noted in the second paragraph of this section, there is a
function $f$ in the class $S$ which maps $\B(1)$ onto the domain
$\Omega = \{z\in \C: 1/z \in K^c\}$. Denoting the Green's function
of $\Omega$ with pole at $0$ by $g(z,0,\Omega)$, the conformal
invariance of Green's functions shows that
$$g(z,0,\Omega) = g(1/z) , \quad z \in \C.$$
Let $\widetilde{\Omega},\widetilde{f}$ and $\widetilde{G}(z, 0,
\widetilde{\Omega})$ be the correponding objects for
$\widetilde{L}$. Then $\widetilde{\Omega} = \C\setminus
[1/4,\infty)$ and $\widetilde{f}(z) = \frac{z}{(1+z)^2}$, the Koebe
function with omitted set on the positive real axis, and
$$g(z,0,\widetilde{\Omega}) = \widetilde{G}(1/z), \quad z \in \C.$$

Thus, the conclusion of Theorem~\ref{logmoment} will hold if for
every $r\in (0,\infty)$,
\begin{equation}
\label{6d} \int_{-\pi}^{\pi}  g(re^{i\theta},0,\Omega)\,d\theta \leq
\int_{-\pi}^{\pi} g(re^{i\theta},0,\widetilde{\Omega})\,d\theta.
\end{equation}

But this inequality is true, since it is the special case $\varphi =
\pi$ in inequality (35) of \cite{Ba1}. (The functions called there
$u^*(re^{i\pi})$ and $v^*(re^{i\pi})$ equal the left and right sides
of \eqref{6d}, respectively.) Theorem~\ref{logmoment} is proved.
\end{proof}

\begin{corollary} \label{symmetric}
Suppose that $K \subset \C$ is compact, connected, contains the
origin and satisfies $\capac (K)=1$, and in addition that $K$ is
symmetric with respect to the origin. Then for every convex function
$\phi:\R \to \R$,  we have
$$\int_K \phi(\lz) \, d\mu_K(z) \leq \int_L \phi(\lz) \, d\mu_L(z).$$
\end{corollary}

Here, as before,  $L = [-2,2]$.

\begin{proof}
To prove the Corollary, use the same construction as in the proofs
of Corollaries~\ref{positivezeros} and \ref{average}. That is, let
$\widetilde{K} = \{z^2:z\in K\}$. Then $\widetilde{K}$ saisfies the
hypotheses of Theorem~\ref{logmoment}, and $\mu_{\widetilde{K}}$ is
the push forward of $\mu_K$ by the map $z \mapsto z^2$. Thus by
Theorem~\ref{logmoment},
\begin{align*}
\int_K \phi(\lz) \, d\mu_K(z)
& = \int_{\widetilde{K}} \phi\left(\frac{1}{2}\lz\right) \, d\mu_{\widetilde{K}}(z) \\
& \leq \int_{\widetilde{L}} \phi\left(\frac{1}{2}\lz\right) \,
d\mu_{\widetilde{L}}(z) = \int_L \phi(\lz) \, d\mu_L(z).
\end{align*}
The inequality in the middle is justified since $\phi(\frac{1}{2} \,
\cdot)$ is convex. Corollary~\ref{symmetric} is proved.
\end{proof}

There are strict inequality statements for Theorem~\ref{logmoment}
and Corollary~\ref{symmetric}, for which we refer to \cite{L}.

To get a closer parallel to Theorem~\ref{Kconn}, it would be nice if
in Corollary~\ref{symmetric} we could drop the symmetry assumption
on $K$ and replace it by the much weaker assumption that the
conformal centroid of $K$ is at the origin. But the example below
shows that no such result can exist.

\begin{example} \rm Hayman \cite[p. 262]{Ha}
built on work of Jenkins \cite{Je} and showed existence of a map
$f(z) = z + \sum_{n=2}^{\infty} A_n z^n$ in the class $S$ for which
$A_2 = 0$ and $M(r,f)\sim c(1-r)^{-2}$ as $r\to 1$, where $M(r,f) =
\max_{\theta}|f(re^{i\theta})|$ and $c$ is some positive constant.
Let $K = \{1/z : z \notin f(\B(1))\}$. Then $K$ is compact and
connected, contains the origin, and $\capac (K)=1$. The Green's
function $g(z) = g(z,\infty, K^c)$ is related to $f$ by $g(z) = \log
1/|f^{-1}(1/z)|$, where $f^{-1}$ is the inverse function of $f$,
from which one calculates that
$$a_1 = -A_2 = 0.$$
We saw in \S\ref{prelims} that $a_1$ is the conformal centroid of
$K$, and thus the conformal centroid of $K$ is $0$.

The behavior of $M(r,f)$ as $r\to 1$ implies that $M(r,g)\sim
c_1r^{1/2}$ as $r\to 0$.  Since $g$ is subharmonic in $\C$, it
follows that $g$ is majorized in any disk by its Poisson integral
over the boundary. Thus,
$$c_2r^{1/2} \leq M(r,g) \leq \frac{3}{2\pi}\int_{-\pi}^{\pi} g(2re^{i\theta})\,d\theta.$$

On the other hand, the Green's function $G$ of $L^c$ satisfies
$M(r,G)\leq c_3 r$ for all $r\in [0,\infty)$. We conclude that
$$\int_{-\pi}^{\pi}  g(re^{i\theta})\,d\theta > \int_{-\pi}^{\pi} G(re^{i\theta})\,d\theta , \quad r \in (0,r_0),$$
for some $r_0 \in (0,1)$.

Take a smooth, convex $\phi$ which is constant on $(-\infty, 2\log
r_0)$, strictly convex on $(2\log r_0, \log r_0)$, and linear on
$(\log r_0, \infty)$. Then Proposition~\ref{representation} gives
$$\int_K \phi(\lz) \, d\mu_K(z) > \int_L \phi(\lz) \, d\mu_L(z),$$
which is the reverse of the moment inequality we might have hoped
would be true.
\end{example}

\vspace{6pt} This example shows the full analogue of
Corollary~\ref{symmetric} does not hold if the symmetry constraint
is relaxed to the centroid constraint. We now propose a substitute,
``averaged'' result. Assume that $K$ is compact, connected, contains
the origin, has $\capac (K)=1$, and also satisfies the centroid
constraint
$$\int_K z \, d\mu_K(z) = 0.$$

Then, as noted in \S\ref{Kconn_proof}, $K\subset \overline{\B(2)}$,
and thus the formula in Proposition~\ref{representation} is valid
for all $R\geq 2$. Fix $R\geq  2$, and define
\begin{align*}
I(r) = I(r,K)
& = \frac{1}{2\pi} \int_{-\pi}^{\pi} g(re^{i\theta})\,d\theta, && r\in [0,\infty), \\
J(r) = J(r,K)
& = \int_{r}^R I(t,K)\,\frac{dt}{t}, && r\in [0,R].
\end{align*}

\begin{conjecture} \label{conj1}
Suppose $K \subset \C$ is compact and connected with $\capac (K)=1$,
and that $0 \in K$ and the conformal centroid of $K$ lies at the
origin. Then for all $R\geq 2$ we have
$$J(r,K)\leq J(r,L), \quad r\in [0,R].$$
\end{conjecture}

An equivalent inequality is
$$\int_{r<|z|<R} g(z)\,|z|^{-2}\,dx\,dy \leq \int_{r<|z|<R} G(z)\,dx\,dy, \quad r\in [0,R].$$
Since $I(r,K) = \log r$ for $r\geq 2$, by \eqref{6a}, it follows
that $I(r,K) = I(r,L)$ for $r\geq 2$ and hence another equivalent
inequality is
$$\int_{r<|z| < \infty} [g(z) - G(z)]|z|^{-2}\,dx\,dy \leq 0, \quad r\in [0,\infty).$$

There is still another equivalent version of Conjecture~\ref{conj1}
involving functions $\phi$, which we'll call Conjecture~\ref{conj2}.

\begin{conjecture} \label{conj2}
Suppose $K \subset \C$ is compact and connected with $\capac (K)=1$,
and that $0 \in K$ and the conformal centroid of $K$ lies at the
origin. Then for all functions $\phi\in C^1(\R)$ such that both
$\phi$ and $\phi^\prime$ are convex, we have
$$\int_K \phi(\lz) \, d\mu_K \leq \int_L \phi(\lz) \, d\mu_L.$$
\end{conjecture}

To see the equivalence, first reduce to the case of smooth $\phi$
with $\phi \equiv 0$ near $-\infty$, by arguing as in the proof of
Theorem~\ref{logmoment}. Then go to Proposition~\ref{representation}
and express the integral over $\B(R)$ in polar coordinates, and
integrate it by parts with respect to $r$. The resulting formula is
\begin{align*}
\int_K \phi(\lz & ) \, d\mu_K(z) \\
& = \int_0^R \phi^{\prime \prime \prime}(\log t)\,\frac{J(t)}{t}\,dt
+ \phi(\log R) - \phi^\prime (\log R) \log R,
\end{align*}
where $R\geq 2$. Now it is immediate that Conjecture~\ref{conj1}
implies Conjecture~\ref{conj2}. As for the converse, one need only
take $\phi(t) = [(t-\log r)^+]^2$, noting $\phi \equiv 0$ near
$-\infty$ and $\phi^{\prime \prime \prime}$ is a positive point mass
at $\log r$.

To conclude, we describe two special cases of Conjecture~\ref{conj2}
which have appeared in the literature as separate conjectures.

The first concerns the class $\Sigma_0$ of all univalent meromorphic
functions $F$ in the exterior $\B^c$ of the unit disk $\B$, with
$F(z) = z + O(z^{-1})$ as $z\to \infty$. The function $F_0(z) =
z+z^{-1}$ belongs to $\Sigma_0$, and maps the exterior of the unit
disk onto the domain $L^c = \overline{\C}\setminus [-2,2]$.

\begin{conjecture}[Pommerenke \cite{Po}]
If $F\in \Sigma_0$ and $0 \in K = F(\B^c)^c$, then
$$\frac{1}{2\pi} \int_{-\pi}^{\pi} |F(e^{i\theta})|\,d\theta \leq
\frac{1}{2\pi} \int_{-\pi}^{\pi} |F_0(e^{i\theta})|\,d\theta =
\frac{4}{\pi}.$$
\end{conjecture}

The best known estimate \cite{Po} is $\frac{1}{2\pi}
\int_{-\pi}^{\pi}|F(e^{i\theta})|\,d\theta \leq 4.02/\pi$. One would
like to replace $4.02$ by $4$.

Note in Pommerenke's Conjecture that $\capac (K$)= 1 and $K$
satisfies the other hypotheses of Conjecture~\ref{conj2}. Moreover,
$d\mu_K$ is the harmonic measure of $K^c=F(\B^c)$ at $\infty$  and
$\frac{d\theta}{2\pi}$ is the harmonic measure of $\B^c$ at $\infty$
(see \cite{Ran}). By conformal invariance of
harmonic measure, we have
\begin{align*}
\frac{1}{2\pi} \int_{-\pi}^{\pi}|F(e^{i\theta})|\,d\theta & = \int_K
|z| \, d\mu_K , \\
\frac{1}{2\pi} \int_{-\pi}^{\pi}|F_0(e^{i\theta})|\,d\theta & =
\int_L |z| \, d\mu_L .
\end{align*}
Thus if Conjecture~\ref{conj2} is true with $\phi(x) = e^x$, then so
is Pommerenke's Conjecture.

Incidentally, the case $\phi(x) = e^{2x}$ of Conjecture~\ref{conj2}
says
\[
\int_K |z|^2 \, d\mu_K \leq \int_L |z|^2 \, d\mu_L ,
\]
which is equivalent as above to
\[
\frac{1}{2\pi} \int_{-\pi}^{\pi} |F(e^{i\theta})|^2 \,d\theta \leq
\frac{1}{2\pi} \int_{-\pi}^{\pi} |F_0(e^{i\theta})|^2 \,d\theta = 2
.
\]
This case of Conjecture~\ref{conj2} can be proved as follows: write
$F(z)=z+\sum_{n=1}^\infty b_n z^{-n}$ and observe
\begin{align*}
\frac{1}{2\pi} \int_{-\pi}^{\pi} |F(e^{i\theta})|^2 \,d\theta
& = 1 + \sum_{n=1}^\infty |b_n|^2 \\
& \leq 1 + \sum_{n=1}^\infty n|b_n|^2 \\
& \leq 2
\end{align*}
by the area theorem \cite[Theorem~1.3]{Pobook} . Clearly equality
holds if and only if $|b_1|=1$ and $b_n=0$ for all $n \geq 2$, which
means $K$ is a rotate of $L$.

The second special case of Conjecture~\ref{conj2} concerns norms of
polynomials. Let $M_K$ be the smallest number $M$ such that
$$\prod_{j=1}^m \|p_j\|_K \leq M^n \|p\|_K$$
for all polynomials $p$ of degree $n\geq 1$ and all polynomials
$p_1,\dots,p_m$ such that $\prod_{j=1}^m p_j = p$. Here
$\|\cdot\|_K$ denotes the sup norm on $K$. The constant $M_K$ was
evaluated in \cite{Pr} as
\[
M_K = \frac{\exp\left(\dis\int_K \log d_K(z) \, d\mu_K
(z)\right)}{\capac (K)}
\]
where $d_K(z) = \max_{t \in K} |z -t|$ is the farthest point
distance function for $K$. (Further properties of $d_K$ have been
studied by Laugesen and Pritsker \cite{LP}, and Gardiner and Netuka
\cite{GN,GN2}.)

The following natural extremal conjecture for $M_K$ was stated in
Pritsker and Ruscheweyh's paper \cite{PR}:

\begin{conjecture} \label{PRconj}
For all compact connected $K\subset \C$ with more than one point, we
have $M_K \leq M_L$.
\end{conjecture}

The constant $M_K$ is invariant under similarity transformations,
and so when studying the conjecture it suffices to assume $\capac
(K) = 1$ and that the conformal centroid of $K$ lies at the origin.

Assuming in addition that $K$ contains its conformal centroid, the
authors of \cite{PR} showed that $M_K < (1.022)M_L$. In the
direction of the conjectured sharp bound (with constant $1$), they
observed
$$\log M_K \leq \int_K \log (2+|z|) \, d\mu_K(z),$$ with
equality when $K = L$.

Now, the function $\phi(t) = \log (2+e^t)$ is convex on $\R$, and
$\phi^\prime$ is convex on $(-\infty, \log 2]$. Replacing $\phi$ on
$(\log 2,\infty)$ by an appropriate quadratic, we obtain a function
$\widetilde{\phi}$ which, along with its derivative, is convex on
all of $\R$.  Suppose Conjecture~\ref{conj2} is true. Then the
inequality in it holds with $\widetilde{\phi}$ in place of $\phi$.
Moreover, $K\subset \overline{\B(2)}$, and so the integrals are the
same for $\phi$ and $\widetilde{\phi}$, which would establish
Conjecture~\ref{PRconj}.

\end{document}